\documentclass[12pt]{article}
\newtheorem{theorem}{Theorem}[section]
\newtheorem{prop}[theorem]{Proposition}
\newtheorem{cor}[theorem]{Corollary}

\newtheorem{lemma}[theorem]{Lemma}
\newtheorem{unnumber}{}

\newtheorem{prepf}[unnumber]{Proof}
\newenvironment{pf}{\prepf\rm}{\endprepf}
\newtheorem{preex}{Example}
\newenvironment{example}{\preex\rm}{\endpreex}

\usepackage{amsmath}
\newcommand{\Aut}{\mathop{\mathrm{Aut}}\nolimits}

\begin{document}
\title{The cycle polynomial of a permutation group}
\author{Peter J. Cameron and Jason Semeraro}
\date{}
\maketitle
\maketitle

\begin{abstract}
The cycle polynomial of a finite permutation group $G$ is the generating
function for the number of elements of $G$ with a given number of cycles:
\[F_G(x) = \sum_{g\in G}x^{c(g)},\]
where $c(g)$ is the number of cycles of $g$ on $\Omega$. In the first part
of the paper, we develop basic properties of this polynomial, and give a
number of examples.

In the 1970s, Richard Stanley introduced the notion of reciprocity for
pairs of combinatorial polynomials. We show that, in a considerable number
of cases, there is a polynomial in the reciprocal relation to the cycle 
polynomial of $G$; this is the \emph{orbital chromatic polynomial} of 
$\Gamma$ and $G$, where $\Gamma$ is a $G$-invariant graph, introduced by the
first author, Jackson and Rudd. We pose the general problem of finding all
such reciprocal pairs, and give a number of examples and characterisations:
the latter include the cases where $\Gamma$ is a complete or null graph or
a tree.

The paper concludes with some comments on other polynomials associated with
a permutation group.
\end{abstract}

\section{The cycle polynomial and its properties}\label{sec:cyclebasic}

Define the \emph{cycle polynomial} of a permutation group $G$ acting on a set
$\Omega$ of size $n$ to be
\[F_G(x) = \sum_{g\in G}x^{c(g)},\]
where $c(g)$ is the number of cycles of $g$ on $\Omega$ (including fixed
points).

Clearly the cycle polynomial is a monic polynomial of degree $n$.

\begin{prop}
If $a$ is an integer, then $F_G(a)$ is a multiple of $|G|$.
\end{prop}

\begin{pf}
Consider the set of colourings of $\Omega$ with $a$ colours (that is, functions
from $\Omega$ to $\{1,\ldots,a\}$. There is a natural action of $G$ on this
set. A colouring is fixed by an element $g\in G$ if and only if it is constant
on the cycles of $g$; so there are $a^{c(g)}$ colourings fixed by $g$. Now the
Orbit-counting Lemma shows that the number of orbits of $G$ on colourings is
\[\frac{1}{|G|}\sum_{g\in G}a^{c(g)};\]
and this number is clearly a positive integer. The fact that $f(a)$ is an integer for all $a$  follows from \cite[Proposition 1.4.2]{rps}.
\end{pf}

\begin{prop}
$F_G(0)=0$; $F_G(1)=|G|$; and $F_G(2)\ge(n+1)|G|$, with equality if and only
if $G$ is set-transitive.
\end{prop}

\begin{pf}
The first assertion is clear.

There is only one colouring with a single colour. 

If there are two colours, say red and blue, then the number of orbits on
colourings is equal to the number of orbits on (red) subsets of $\Omega$.
There are $n+1$ possible cardinalities of subsets, and so at least $n+1$
orbits, with equality if and only if $G$ is transitive on sets of size $i$
for $0\le i\le n$. (The set-transitive groups were determined by Beaumont and
Petersen \cite{bp}; there are only the symmetric and alternating groups and
four others with $n=5,6,9,9$.)
\end{pf}

Now we consider values of $F_G$ on negative integers. Note that the \emph{sign}
of the permutation $g$ is $(-1)^{n-c(g)}$; a permutation is even or odd 
according as its sign is $+1$ or $-1$. If $G$ contains odd permutations, then
the even permutations in $G$ form a subgroup of index~$2$.

\begin{prop}
If $G$ contains no odd permutations, then $F_G$ is an even or odd function 
according as $n$ is even or odd; in other words,
\[F_G(-x)=(-1)^nF_G(x).\]
\label{p:parity}
\end{prop}

\begin{pf}
The degrees of all terms in $F_G$ are congruent to $n$ mod~$2$.
\end{pf}

In particular, we see that if $G$ contains no odd permutations, then $F_G(x)$
vanishes only at $x=0$. However, for permutation groups containing odd
permutations, there may be negative roots of $F_G$.

\begin{theorem}
Suppose that $G$ contains odd permutations, and let $N$ be the subgroup of
even permutations in $G$. Then, for any positive integer $a$,
we have $0\le(-1)^nF_G(-a)<F_G(a)$, with equality if and only if $G$ and $N$
have the same number of orbits on colourings of $\Omega$ with $a$ colours.
\end{theorem}

\begin{pf}
We have
\begin{eqnarray*}
\sum_{g\in N}a^{c(g)} + \sum_{g\in G\setminus N}a^{c(g)} &=& F_G(a),\\
\sum_{g\in N}a^{c(g)} - \sum_{g\in G\setminus N}a^{c(g)} &=& (-1)^nF_G(-a).
\end{eqnarray*}
Moreover, $N$ has at least as many orbits on colourings as does $G$, and so
we have
\[\sum_{g\in N}a^{c(g)}/|N| \ge \sum_{g\in G}a^{c(g)}/|G|.\]
So $F_N(a)$ is at least half of $F_G(a)$, showing that the left-hand side
of the second equation above is at least zero. This proves the inequality;
we see that equality holds if and only if $G$ and $N$ have equally many orbits
on colourings, as required.
\end{pf}

\begin{prop}
If $G$ is a permutation group containing odd permutations, then the set of
negative integer roots of $F_G$ consists of all integers $\{-1,-2,\ldots,-a\}$
for some $a\ge1$.
\end{prop}

\begin{pf}
$F_G(-1)=0$, since $G$ and $N$ have equally many orbits (namely $1$) on
colourings with a single colour.

Now suppose that $F_G(-a)=0$, so that $G$ and $N$ have equally many orbits on
colourings with $a$ colours; thus every $G$-orbit is an $N$-orbit. Now every
colouring with $a-1$ colours is a colouring with $a$ colours, in which the
last colour is not used; so every $G$-orbit on colourings with $a-1$ colours
is an $N$-orbit, and so $F_G(-a+1)=0$. The result follows.
\end{pf}

The property of having a root $-a$ is preserved by overgroups:

\begin{prop}\label{p:overgroups}
Suppose that $G_1$ and $G_2$ are permutation groups on the same set, with
$G_1\le G_2$. Suppose that $F_{G_1}(-a)=0$, for some positive integer $a$.
Then also $F_{G_2}(-a)=0$.
\end{prop}

\begin{pf}
It follows from the assumption that $G_1$ (and hence also $G_2$) contains odd
permutations. Let $N_1$ and $N_2$ be the subgroups of even permutations in
$G_1$ and $G_2$ respectively. Then $N_2\cap G_1=N_1$, and so $N_2G_1=G_2$.
By assumption, $G_1$ and $N_1$ have the same orbits on $a$-colourings. Let
$K$ be an $a$-colouring, and $g\in G_2$; write $g=hg'$, with $h\in N_2$ and
$g'\in G_1$. Now $Kh$ and $Khg'$ are in the same $G_1$-orbit, and hence in the
same $N_1$-orbit; so there exists $h\in N_1$ with $Kg=Khg'=Khh'$. Since
$hh'\in N_2$, we see that the $G_2$-orbits and $N_2$-orbits on $a$-colourings
are the same. Hence $F_{G_2}(-a)=0$.
\end{pf}

The cycle polynomial has nice behaviour under direct product, which shows that
the product of having negative integer roots is preserved by direct product.

\begin{prop}\label{p:prod}
Let $G_1$ and $G_2$ be permutation groups on disjoint sets $\Omega_1$ and
$\Omega_2$. Let $G=G_1\times G_2$ acting on $\Omega_1\cup\Omega_2$. Then
\[F_G(x)=F_{G_1}(x)\cdot F_{G_2}(x).\]
In particular, the set of roots of $F_G$ is the union of the sets of roots of
$F_{G_1}$ and $F_{G_2}$.
\end{prop}

\begin{pf}
This can be done by a calculation, but here is a more conceptual proof. It
suffices to prove the result when a positive integer $a$ is substituted for
$x$. Now a $G$-orbit on $a$-colourings is obtained by combining a $G_1$-orbit
on colourings of $\Omega_1$ with a $G_2$-orbit of colourings of $\Omega_2$; so
the number of orbits is the product of the numbers for $G_1$ and $G_2$.
\end{pf}

The result for the wreath product, in its imprimitive action, is obtained in a
similar way.

\begin{prop}\label{p:wreath}
\[F_{G\wr H}(x) = |G|^mF_H(F_G(x)/|G|).\]
\end{prop}

\begin{pf} Again it suffices to prove that, for any positive integer $a$, the
equation is valid with $a$ substituted for $x$.

Let $\Delta$ be the domain of $H$, with $|\Delta|=m$. An orbit of the base
group $G^m$ on $a$-colourings is an $m$-tuple of $G$-orbits on $a$-colourings,
which we can regard as a colouring of $\Delta$, from a set of colours whose
cardinality is the number $F_G(a)/|G|$ of $G$-orbits on $a$-colourings. Then an
orbit of the wreath product on $a$-colourings is given by an orbit of $H$ on
these $F_G(a)/|G|$-colourings, and so the number of orbits is
$(1/|H|)F_H(F_G(a)/|G|)$. Multiplying by $|G\wr H|=|G|^m|H|$ gives the
result.
\end{pf}

\begin{cor}
If $n$ is odd, $m > 1$, and $G=S_n\wr S_m$, then $F_G(x)$ has roots $-1,\ldots,-n$.
\end{cor}

\begin{pf}
Clearly $F_{S_n}(x)=x(x+1)\cdots(x+n-1)$ divides $F_G(x)$, so we have roots
$-1,\ldots,-n+1$. Also, there is a factor
\[|S_n|(F_{S_n}(x)/|S_n|+1) = x(x+1)\cdots(x-n+1)+n!.\]
Substituting $x=-n$ and recalling that $n$ is odd, this is $-n!+n!=0$.
\end{pf}

The next corollary shows that there are imprimitive groups with arbitrarily
large negative roots.
 
\begin{cor}
Let $n$ be odd and let $G$ be a permutation group of degree $n$ which contains
no odd permutations. Suppose that $F_G(a)/|G|=k$. Then, for $m>k$, the
polynomial $F_{G\wr S_m}(x)$ has a root $-a$.
\end{cor}

\begin{pf}
By Proposition~\ref{p:parity}, $F_G(-a)/|G|=-k$. Now for $m>k$, the
polynomial $F_{S_m}(x)$ has a factor $x+k$; the expression for $F_{G\wr S_m}$
shows that this polynomial vanishes when $x=-a$.
\end{pf}

The main question which has not been investigated here is:
\begin{quote}
What about non-integer roots?
\end{quote}

It is clear that $F_G(x)$ has no positive real roots; so if $G$ contains no
odd permutations, then $F_G(x)$ has no real roots at all, by
Proposition~\ref{p:parity}. Is there any restriction on where the non-real
roots lie? For example computer calculations suggest that $F_{A_n}(a)=0$ implies $\operatorname{Re}(a)=0$, but we have been unable to prove this.

\section{Some examples}

\begin{prop}\label{p:sym}
For each $n$ we have, $$F_{S_n}(x)=\prod_{i=0}^{n-1} (x+i).$$
\end{prop}

\begin{pf}
By induction and Proposition \ref{p:overgroups}, $F_{S_n}(-a)=0$ for each $0 \le a \le n-2$. Since $F_{S_n}(x)$ is a polynomial of degree $n$ and the coefficient of $x$ in $F_{S_n}(x)$ is $(n-1)!$ we must have that $n-1$ is the remaining root.

Several other proofs of this result are possible.
We can observe that the number of orbits of $S_n$ on
$a$-colourings of $\{1,\ldots,n\}$ (equivalently, $n$-tuples chosen from the
set of colours with order unimportant and repetitions allowed) is
$n+a-1\choose n$. Or we can use the fact that the number of permutations of
$\{1,\ldots,n\}$ with $k$ cycles is the unsigned \emph{Stirling number of the
first kind} $u(n,k)$, whose generating function is well known to be
\[\sum_{k=1}^nu(n,k)x^k=x(x+1)\cdots(x+n-1).\]
\end{pf}

\begin{prop}
For each $n$ we have $$F_{C_n}(x)=\sum_{d \mid n} \phi(d)x^{n/d},$$ where $\phi$ is  Euler's totient function.
\end{prop}

\begin{pf}
Clear.
\end{pf}

\begin{prop}\label{prop:fg}
Let $p$ be an odd prime and $G$ be the group $PGL_2(p)$ acting on a set of size $p+1$. Then $F_G(x)$ is given by
$$\frac{p(p-1)}{2} \cdot F_{C_{p+1}}(x) +  \frac{p(p+1)}{2} \cdot x^2F_{C_{p-1}}(x) + (p^2-1)(x^2-x^{p+1})$$
\end{prop}

\begin{pf}
We count the elements in $G$ of order $k$. When $k=2$, there are two conjugacy classes of involutions of lengths $\frac{p(p-1)}{2}$ and $\frac{p(p+1)}{2}$ with the first class consisting of elements which act fixed-point freely and the second consisting of elements which fix exactly two points. For $\epsilon = \pm 1$, and each non-trivial odd divisor $k$ of $p+\epsilon$, there are $$\phi(k)\frac{p(p-\epsilon)}{2}$$ elements of order $k$ with 
$$
c(g)=
\begin{cases}
 \frac{p+1}{k},& \mbox{ if $\epsilon = 1$ }\\
 \frac{p-1}{k} + 2,& \mbox{ if $\epsilon = -1$ }.\\
\end{cases}
$$
Finally, there is a single conjugacy class of elements of order $p$ of length $p^2-1$. Putting this altogether, we obtain
\begin{eqnarray*}
F_G(x) &=& x^{p+1}+\frac{p(p-1)}{2}x^{\frac{p+1}{2}}+\frac{p(p+1)}{2}x^{\frac{p+3}{2}}+ \frac{p(p-1)}{2} \sum_{1,2 \neq k \mid p+1} \phi(k)x^{\frac{p+1}{k}}\\
&&\quad+\frac{p(p+1)}{2} \sum_{1,2 \neq k \mid p-1} \phi(k)x^{\frac{p-1}{k}+2}+(p^2-1)x^2.
\end{eqnarray*}
Simplifying gives the result.
\end{pf}

\section{Reciprocal pairs}

Richard Stanley, in a 1974 paper \cite{stanley}, explained (polynomial)
combinatorial reciprocity thus:
\begin{quote}
A \emph{polynomial reciprocity theorem} takes the following form. Two
combinatorially defined sequences $S_1$, $S_2$, \ldots\ and $\bar S_1$,
$\bar S_2$, \ldots\ of finite sets are given, so that the functions
$f(n)=|S_n|$ and $\bar f(n)=|\bar S_n|$ are polynomials in $n$ for all
integers $n\ge1$. One then concludes that $\bar f(n)=(-1)^df(-n)$, where
$d=\deg f$.
\end{quote}

We will see that, in a number of cases, the cycle polynomial satisfies a
reciprocity theorem.

\subsection{The orbital chromatic polynomial}

First, we define the polynomial which will serve as the reciprocal polynomial
in these cases. A \emph{(proper) colouring} of a graph $\Gamma$ with $q$
colours is a map from the vertices of $\Gamma$ to the set of colours having
the property that adjacent vertices receive different colours. Note that,
if $\Gamma$ contains a loop (an edge joining a vertex to itself), then it has
no proper colourings. Birkhoff observed that, if there are no loops, then
the number of colourings with $q$ colours is the evaluation at $q$ of a
monic polynomial $P_\Gamma(x)$ of degree equal to the number of vertices, the
\emph{chromatic polynomial} of the graph.

Now suppose that $G$ is a group of automorphisms of $\Gamma$. For $g\in G$,
let $\Gamma/g$ denote the graph obtained by ``contracting'' each cycle of $g$
to a single vertex; two vertices are joined by an edge if there is an edge
of $\Gamma$ joining vertices in the corresponding cycles. The chromatic
polynomial $P_{\Gamma/G}(q)$ counts proper $q$-colourings of $\Gamma$ fixed
by $g$. If any cycle of $g$ contains an edge, then $\Gamma/g$ has a loop, and
$P_{\Gamma/g}=0$. Now (with a small modification of the definition
in~\cite{cjr}) we define the \emph{orbital chromatic polynomial} of the pair
$(\Gamma,G)$ to be
\[P_{\Gamma,G}(x)=\sum_{g\in G}P_{\Gamma/g}(x).\]
The Orbit-Counting Lemma immediately shows that $P_{\Gamma,G}(q)/|G|$ is equal to
the number of $G$-orbits on proper $q$-colourings of $\Gamma$.

Now, motivated by Stanley's definition, we say that the pair $(\Gamma,G)$,
where $\Gamma$ is a graph and $G$ a group of automorphisms of $\Gamma$, is
a \emph{reciprocal pair} if
\[P_{\Gamma,G}(x)=(-1)^nF_G(-x),\]
where $n$ is the number of vertices of $\Gamma$.

\paragraph{Problem} Find all reciprocal pairs.

\medskip

This problem is interesting because, as we will see, there are a substantial
number of such pairs, for reasons not fully understood. In the remainder of
the paper,  we present the evidence for this, and some preliminary results
on the above problem.

A basic result about reciprocal pairs is the following.

\begin{lemma}\label{lem:sumtrans}
Suppose that $(G,\Gamma)$ is a reciprocal pair. Then the number of edges of
$\Gamma$ is the sum of the number of transpositions in $G$ and the number of
transpositions $(i,j)$ in $G$ for which $\{i,j\}$ is a non-edge.
\end{lemma}

\begin{pf}
Whitney~\cite{whitney} showed that the leading terms in the chromatic
polynomial of a graph $\Gamma$ with $n$ vertices and $m$ edges are
$x^n-mx^{n-1}$. The only additional contributions to the coefficient of
$x^{n-1}$ in $P_{\Gamma,G}(x)$
come from elements of $G$ with $n-1$ cycles, that is, transpositions. A
transposition $(i,j)$ makes a non-zero contribution if and only if $\{i,j\}$
is a non-edge. So the coefficient of $x^{n-1}$ is $-m+t^0(G)$, where $t^0(G)$
is the number of transpositions with this property.

On the other hand, the coefficient of $x^{n-1}$ in $F_G(x)$ is the number of
permutations in $G$ with $n-1$ cycles, that is, the total number $t(G)$ of
transpositions. So the coefficient in $(-1)^nF_G(-x)$ is $-t(G)$.

Equating the two expressions gives $m=t(G)+t^0(G)$, as required.
\end{pf}

We remark that the converse to Lemma \ref{lem:sumtrans} does not hold: consider the group $G=S_3 \wr S_3$ acting on 3 copies of $K_3$ (see
Proposition~\ref{p:wreath}: but note that $(3K_3, S_3\wr C_3)$ is a reciprocal
pair, by Proposition~\ref{p:wreath_recip}). We also observe the following corollary to Lemma \ref{lem:sumtrans}.

\begin{cor}
If $\Gamma$ is not a complete graph and $(\Gamma,G)$ is a reciprocal pair then $\Gamma$ has at most $\frac{(n-1)^2}{2}$ edges.
\end{cor}

\begin{pf}
If $\Gamma$ has ${n \choose 2}-\delta$ edges then by Lemma \ref{lem:sumtrans}, $${n \choose 2}-\delta = t(G)+t^0(G) \le t(G) + \delta.$$ If $0  < \delta < \frac{n-1}{2}$ then $$t(G) \ge {n \choose 2}-2\delta > {n \choose 2} - (n-1) ={n-1 \choose 2}.$$ It is well-known that a permutation group of degree $n$ containing at least ${n-1 \choose 2}+1$ transpositions must be the full symmetric group. But this implies that $\Gamma$ is a complete graph, a contradiction. 
\end{pf}

\medskip

According to Lemma \ref{lem:sumtrans}, if $\Gamma$ is not a null graph and
$(\Gamma,G)$ is a reciprocal pair, then $G$ contains transpositions. Now
as is well known, if a subgroup $G$ of $S_n$ contains a transposition, then
the transpositions generate a normal subgroup $N$ which is the direct product
of symmetric groups whose degrees sum to $n$. (Some degrees may be $1$, in
which case the corresponding factor is absent.) 

A $G$-invariant graph must induce a complete or null graph on each of these
sets. Moreover, between any two such sets, we have either all possible edges
or no edges.

Suppose that $n_1,\ldots,n_r$ are the $N$-orbits carrying complete graphs and
$m_1,\ldots,m_s$ the orbits containing null graphs, then
Lemma~\ref{lem:sumtrans} shows that the total number of edges of the graph
is
\[\sum_{i=1}^r{n_i\choose 2}+2\sum_{j=1}^s{m_j\choose 2}.\]
The first term counts edges within $N$-orbits, so the second term counts
edges between different $N$-orbits.

\subsection{Examples}

\begin{prop}\label{p:ex1}
\begin{itemize}
\item[(a)] Let $\Gamma$ be a null graph, and $G$ a subgroup of the symmetric group
$S_n$. Then $P_{\Gamma,G}(x)=F_G(x)$.
\item[(b)] Let $\Gamma$ be a complete graph, and $G$ a subgroup of the symmetric
group $S_n$. Then $P_{\Gamma,G}(x)=x(x+1)\cdots(x+n-1)$, independent of
$\Gamma$.
\end{itemize}
\end{prop}

\begin{pf}
(a) The chromatic polynomial of a null graph on $n$ vertices is $x^n$. So, if
$g\in G$ has $c(g)$ cycles, then $\Gamma/g$ is a null graph on $c(g)$ vertices.
Thus
\[P_{\Gamma,G}(x)=\sum_{g\in G}x^{c(g)}=F_G(x).\]

\medskip

(b) Suppose that $\Gamma$ is a complete graph. Any proper colouring of $\Gamma$
has all colours disjoint, so any permutation group $G$ on proper $a$-colourings
is semiregular, and has $a(a-1)\cdots(a-n+1)/|G|$ orbits. Thus the orbital
chromatic polynomial of $G$ is $x(x-1)\cdots(x-n+1)$, independent of $G$.
\end{pf}

\begin{cor}
\begin{itemize}
\item[(a)] If $\Gamma$ is a null graph, then $(\Gamma,G)$ is a reciprocal pair
if and only if $G$ contains no odd permutations.
\item[(b)] If $\Gamma$ is a complete graph, then $(\Gamma,G)$ is a reciprocal
pair if and only if $G$ is the symmetric group.
\end{itemize}
\end{cor}

\begin{pf}
(a) This follows from Proposition~\ref{p:parity}.

\medskip

(b) We saw in the preceding section that, if $G=S_n$, then
$F_G(x)=x(x+1)\cdots(x+n-1)$. Thus, we see that $(G,\Gamma)$ is a reciprocal
pair if and only if $G$ is the symmetric group.
\end{pf}

\begin{prop}\label{p:direct_recip}
Let $\Gamma$ be the disjoint union of graphs $\Gamma_1,\ldots,\Gamma_r$,
and $G$ the direct product of groups $G_1,\ldots,G_r$, where
$G_i\le\Aut(\Gamma_i)$. Then
\[P_{\Gamma,G}(x)=\prod_{i=1}^rP_{\Gamma_i,G_i}(x).\]
In particular, if $(\Gamma_i,G_i)$ is a reciprocal pair for $i=1,\ldots,r$,
then $(\Gamma,G)$ is a reciprocal pair.
\end{prop}

The proof is straightforward; the last statement follows from
Proposition~\ref{p:prod}. The result for wreath products is similar:

\begin{prop}\label{p:wreath_recip}
Let $\Gamma$ be the disjoint union of $m$ copies of the $n$-vertex graph
$\Delta$. Let $G\le\Aut(\Delta)$, and $H$ a group of permutations of degree
$m$. Then
\[P_{\Gamma,G\wr H}(x) = |G|^mF_H(P_{\Delta,G}(x)/|G|).\]
In particular, if $(\Delta,G)$ is a reciprocal pair and $H$ contains no odd
permutations, then $(\Gamma,G\wr H)$ is a reciprocal pair.
\end{prop}

\begin{pf}
Given $q$ colours, there are $P_{\Delta,G}(q)/|G|$ orbits on colourings of each
copy of $\Delta$; so the overall number of orbits is the same as the number of
orbits of $H$ on an $m$-vertex null graph with $P_{\Delta,G}(q)/|G|$ colours
available.

For the last part, the hypotheses imply that $P_{\Delta,G}(q)=(-1)^nF_G(-q)$,
and that the degree of each term in $F_H$ is congruent to $m$ mod~$2$. So
the expression evaluates to $|G|^mF_H(F_G(-x)/|G|)$ if either $m$ or $n$ is
even, and the negative of this if both are odd; that is,
$(-1)^{mn}|G|^mF_H(F_G(-x)/|G|)$.
\end{pf}

We now give some more examples of reciprocal pairs.

\begin{example}
Let $\Gamma$ be a $4$-cycle, and $G$ its automorphism group, the dihedral
group of order $8$. There are $4$ edges in $\Gamma$, and $2$ transpositions in
$G$, each of which interchanges two non-adjacent points (an opposite pair of
vertices of the $4$-cycle); so the equality of the lemma holds. Direct
calculation shows that
\[P_{\Gamma,G}(x)=x(x-1)(x^2-x+2),\qquad F_G(x)=x(x+1)(x^2+x+2),\]
so $P_{\Gamma,G}(x)=(-1)^nF_G(-x)$ holds in this case.

The $4$-cycle is also the complete bipartite graph $K_{2,2}$. We note that,
for $n>2$, $(K_{n,n},S_n\wr S_2)$ is not a reciprocal pair. This can be seen
from the fact that $F_{S_n\wr S_2}(x)$ has factors $x$, $x+1,\ldots,x+n-1$, whereas $K_{n,n}$ has chromatic number $2$ and so $x-2$ is not a
factor of $P_{K_{n,n},G}(x)$ for any $G\le S_n\wr S_2$. Also  $(K_{n,m},S_n \times S_m)$ is not a reciprocal pair for $n,m > 2$ since the equality $$mn=2\left({m \choose 2}+{n \choose 2}\right)=m(m-1)+n(n-1)$$ of Lemma \ref{lem:sumtrans} is easily seen to have no solutions in this range.

We do not know whether other complete multipartite graphs support
reciprocal pairs.
\end{example}

\begin{example}
Let $\Gamma$ be a path with 3 vertices and $G$ its automorphism group which is cyclic of order 2. Direct calculation shows that 
\[P_{\Gamma,G}(x)=x^2(x-1),\qquad F_G(x)=x^2(x+1).\]
This graph is an example of a star graph, for which we give a complete
analysis in the next section.
\end{example}

\begin{example} Let $\Gamma$ be the disjoint union of $K_m$ and $N_n$ together with all edges in between and set $G=S_m \times S_n \le$ Aut$(\Gamma)$. Then $\Gamma$ has ${m \choose 2} +mn$ edges and $G$ has ${m \choose 2}+{n \choose 2}$ transpositions of which ${n \choose 2}$ correspond to non-edges in $\Gamma$. Thus, according to Lemma \ref{lem:sumtrans} we need $n=m+1$. Now the only elements $g \in G$ which give non-zero contribution to $P_{\Gamma,G}(x)$ lie in the $S_n$ component. We get: $$P_{\Gamma,G}(x)=x(x-1)\cdots (x-(m-1)) \cdot \sum_{g \in S_n} (x-m)^{c(g)}.$$ Using Proposition \ref{p:ex1}(b) and $m=n-1$ this becomes $$(-1)^m F_{S_m}(-x) F_{S_n}(-x) \cdot (-1)^n,$$ which is equal to $-F_G(-x)$ by Proposition \ref{p:prod}.
\end{example}

\section{Reciprocal pairs containing a tree}

In this section we show that the only trees that can occur in a reciprocal
pair are stars, and we determine the groups that can be paired with them. 

\begin{theorem}\label{t:recpairstree}
Suppose $\Gamma$ is a tree and $(\Gamma,G)$ is a reciprocal pair. Let $n$ be the number of vertices in $\Gamma$ and assume $n \ge 3$. The following hold:
\begin{itemize}\itemsep0pt
\item[(a)] $n$ is odd;
\item[(b)] $\Gamma$ is a star;
\item[(c)] $(C_2)^k \le G \le C_2 \wr S_k$ where $n=2k+1$.
\end{itemize}
Conversely any pair $(\Gamma,G)$ which satisfies conditions (a)-(c) is reciprocal.
\end{theorem} 

In what follows we assume the following:
\begin{itemize}\itemsep0pt
\item $\Gamma$ is a tree;
\item $(\Gamma,G)$ is a reciprocal pair;
\item $n$ is the number of vertices in $\Gamma$ and $n \ge 3$.
\end{itemize}

Note that any two vertices interchanged by a transposition are non-adjacent.
For suppose that a transposition flips an edge $\{v,w\}$. If the tree is
central, then there are paths of the same length from the centre to $v$ and
$w$, creating a cycle. If it is bicentral, then the same argument applies
unless $\{v,w\}$ is the central edge, in which case the tree has only two 
vertices, a contradiction.

\begin{lemma}\label{lem:treeequal}
If $n$ is odd then $(\Gamma,G)$ is a reciprocal pair if and only if
\begin{equation}\label{e:eq1} x(F_G(x-1)+F_G(-x))=F_G(-x).
\end{equation}
\end{lemma}

\begin{pf}
 The chromatic polynomial of a tree with $r$ vertices is easily seen to be $x(x-1)^{r-1}$. Hence $$P(\Gamma/g)=x(x-1)^{c(g)-1}$$ for each $g \in G$ and then
 $$P_{\Gamma,G}(x)=\sum_{g \in G} P(\Gamma/g) = \sum_{g \in G} x(x-1)^{c(g)-1} = \frac{x}{x-1}F_G(x-1).$$ Rearranging (and using that $n$ is odd) yields the Lemma.
\end{pf}

\begin{lemma}\label{lem:nodd}
$G$ has $\frac{n-1}{2}$ transpositions; in particular, $n$ is odd.
\end{lemma}

\begin{pf}
As in Lemma \ref{lem:sumtrans}, let $t(G)$ be the number of transpositions in $G$ and $t^0(G)$ be the number of transpositions $(i,j)$ in $G$ for which $i \not\sim j$ in $\Gamma$. If $G$ fixes an edge $(u,v)$ of $\Gamma$ then $(u,v) \in G$ implies $n = 2$, a contradiction. Thus $(u,v) \notin G$, and every transposition in $G$ is a non-edge. Hence $t^0(G)=t(G)$ so $2t(G)=n-1$ by Lemma \ref{lem:sumtrans}.

\end{pf}

\begin{lemma}\label{lem:groupstruct}
$(C_2)^k \le G \le C_2 \wr S_k$ where $n=2k+1$.
\end{lemma}

\begin{pf}
The transpositions in a permutation group $G$ generate a normal subgroup $H$
which is a direct product of symmetric groups. If there are two non-disjoint
transpositions in $G$, one of the direct factors is a symmetric group with
degree at least $3$, and hence $F_H(x)$ has a root $-2$ by Propositions
\ref{p:sym} and \ref{p:prod}. Then by Proposition~\ref{p:overgroups},
$F_G(x)$ has a root $-2$. By Lemma \ref{lem:treeequal} with $x=2$,

$$0=F_G(-2)=2(F_G(1)+F_G(-2))=2F_G(1) = 2|G|,$$ a contradiction. So the transpositions are pairwise disjoint, and generate a subgroup $(C_2)^k$ with $n=2k+1$ by Lemma \ref{lem:nodd}. Thus the conclusion of the lemma holds.
\end{pf}

\begin{lemma}\label{lem:star}
$\Gamma$ is a star;
\end{lemma}

\begin{pf}
Let $v$ be the unique fixed point of $G$. By Lemma \ref{lem:groupstruct}, for each $u \neq v$ there exists a unique vertex $u'$ with $(u,u') \in G$. This is possible only if each $u$ has distance 1 from $v$. Hence $\Gamma$ is a star.
\end{pf}

\begin{pf}[Proof of Theorem \ref{t:recpairstree}]
(a),(b) and (c) follow from Lemmas \ref{lem:nodd}, \ref{lem:star} and \ref{lem:groupstruct} respectively. Conversely, suppose that (a),(b) and (c) hold. Then $G=C_2\wr K$ for some permutation
group $K$ of degree $k$. By Proposition~\ref{p:wreath}, 
$F_G(x)=x \cdot 2^kF_K(x(x+1)/2)$. Now it is clear that $$-\frac{F_G(-x)}{x}=2^k \cdot F_K\left(\frac{x(x-1)}{2}\right)=\frac{F_G(x-1)}{x-1},$$ so that
(\ref{e:eq1}) holds and we deduce from Lemma \ref{lem:treeequal} that $(\Gamma,G)$ is a reciprocal pair. Our proof is complete.
\end{pf}

Given a set of reciprocal pairs $(\Gamma_1,G_1),\ldots (\Gamma_m,G_m)$ with with each $\Gamma_i$ a star we can take direct products and wreath products (using Propositions \ref{p:direct_recip} and \ref{p:wreath_recip}) to obtain reciprocal pairs $(\Gamma,G)$ with $\Gamma$ a forest of stars. We do not know whether all such pairs arise in this way.

\section{Connection with other polynomials}

The \emph{cycle index} of a permutation group $G$ of degree $n$ is a polynomial
in $n$ variables which keeps track of all the cycle lengths of elements, not
just the total number of cycles. If the variables are $s_1,\ldots,s_n$, then
the cycle index is given by
\[Z_G(s_1,\ldots,s_n) = \sum_{g\in G}\prod_{i=1}^n s_i^{c_i(g)},\]
where $c_i(g)$ is the number of cycles of length $i$ in the cycle decomposition
of $G$. (It is customary to divide this polynomial by $|G|$; but, for
compatibility with our earlier polynomial and comparison with other polynomials
we do not do so here.)

Clearly the cycle polynomial is given by $F_G(x)=Z_G(x,x,\ldots,x)$.

Harden and Penman studied the \emph{fixed point polynomial} $P_G(x)$ of a
permutation group $G$ (the generating function for fixed points of elements of
$G$), given by $P_G(x)=Z_G(x,1,\ldots,1)$.

Another related construction is the \emph{Parker vector} of $G$, usually 
presented as an $n$-tuple rather than a polynomial: its $k$th entry is 
obtained by differentiating $Z_G$ with respect to $s_k$, and putting all
variables equal to $1$ (and dividing by $|G|$). It is relevant to computing
the Galois group of an integer polynomial. See~\cite[Section 2.8]{permgps}.

In \cite{cjr}, the orbital chromatic polynomial is extended to an \emph{orbital Tutte polynomial}, for graphs or (more generally) representable matroids. It is not clear whether this is related to the cycle polynomial.

Also, the cycle index of finite permutation groups can be extended to \emph{oligomorphic} (infinite) permutation groups, see \cite[Section 5.7]{permgps}. However, the specialisation which gives the cycle polynomial fails to be defined in the infinite case.

\end{document}